\documentclass[runningheads,a4paper]{llncs}

\usepackage{amssymb}
\usepackage{graphicx}
\usepackage{comment}

\usepackage{url}

\begin{document}

\mainmatter  

\title{A dynamic symbolic geometry environment based on the Gr\"obnerCover algorithm for the computation of geometric loci and envelopes}

\titlerunning{A geometric environment: loci and envelopes}

%
%
\author{Miguel A. Ab\'anades \and Francisco Botana \thanks{Both authors partially supported by the project MTM2011-25816-C02-(01,02) funded by the Spanish \textit{Ministerio de Econom\'ia y Competitividad} and the European Regional Development Fund (ERDF)\dots The final publication is available at http://link.springer.com.}}
\authorrunning{Miguel A. Ab\'anades \and Francisco Botana}
\institute{CES Felipe II, Universidad Complutense de Madrid\\
C/ Capit\'an 39, 28300 Aranjuez, Spain\\
abanades@ajz.ucm.es\\
Depto. de Matem\'atica Aplicada I, Universidad de Vigo\\
Campus A Xunqueira, 36005 Pontevedra, Spain\\
fbotana@uvigo.es}

%

\maketitle

\begin{abstract}
An enhancement of the dynamic geometry system \mbox{GeoGebra} for the automatic symbolic computation of algebraic loci and envelopes is presented. Given a \mbox{GeoGebra} construction, the prototype, after rewriting the construction as a polynomial system in terms of variables and parameters, uses an implementation of the recent Gr\"obnerCover algorithm to obtain the algebraic description of the sought locus/envelope as a locally closed set. The prototype shows the applicability of these techniques in general purpose dynamic geometry systems.
\keywords{Dynamic Geometry, Locus, Envelope, Gr\"obnerCover Algorithm, GeoGebra, Sage}
\end{abstract}

\section{Introduction}

Most dynamic geometry systems (DGS) implement loci generation just from a graphic point of view, returning a locus as a set of points in the screen with no algebraic information. A simple algorithm based on elimination theory to obtain the equation of an algebraic plane curve from its description as a locus set was described in \cite{BotanaValcarceMatCom2003}. This new information expands the algebraic knowledge of the system, allowing further transformations of the construction elements, such as constructing a point on a locus, intersecting the locus with other elements, etc. The same consideration can be made with respect to the envelope of a family of curves. This algebraic approach is a significant improvement over the numeric-graphic method mentioned above. An implementation of the algorithm in a system embedding GeoGebra in the Sage notebook was described at CICM 2011 \cite{BotanaCICM2011}. In fact, the algorithm is already behind the \textit{LocusEquation} command in the beta version of the next version of the DGS GeoGebra \cite{GeoGebra} (see \url{http://wiki.geogebra.org/en/LocusEquation_Command}). It has also recently been implemented by the DGS JSXGraph \cite{jsxgraph} to determine the equation of a locus set using remote computations on a server \cite{JSXGraphADG2010}, an idea previously developed by the authors in \cite{LAD}.  

Unfortunately, this algorithm does not discriminate between regular and special components of a locus (following the definitions in \cite{SendraSendra2000})\footnote{A special component of a locus is basically a one-dimensional subset of the locus corresponding to a single position of the moving point.}. More concretely, the obtained algebraic set may contain extra components sometimes due to the fact that the method returns only Zariski closed sets (i.e. zero sets of polynomials) and sometimes due to degenerate positions in the construction (e.g. two vertices being coincidental for a triangle construction).

There is little that can be done to solve these problems with the simple elimination approach. Concerning degeneration, there is no alternative except explicitly requesting information from the user about the positions producing special components. However, the recent Gr\"obnerCover algorithm \cite{MontesWibmer2010} has opened new possibilities for the automated processing of these problems. From the canonical decomposition of a polynomial system with parameters returned by the algorithm, and following a remark by Tom\'as Recio concerning the dimensions of the spaces of variables and parameters, a protocol has been established to distinguish between regular and special components of a locus set. For example, a circle, a variety of dimension 1, is declared to be a special component of a locus by the protocol if it corresponds to a point, a variety of dimension 0. This heuristic in the protocol improves the automatic determination of loci but does not fully resolve it. It is not difficult to find examples where this general rule does not suit the user's interests. This is a delicate issue because, in some situations, these special components are the relevant parts of the sought set (the study of bisector curves is a good source for such examples). 

As an illustration, let us consider the following problem included in \cite{Guzman2002} together with a remark about its difficult synthetic treatment: {\it Given a triangle $ABC$. Take a point $M$ on $BC$. Consider the orthogonal projections $N$ of $M$ onto $AC$, and $P$ onto $AB$ respectively. The lines $AM$ and $PN$ meet at $X$. What is the locus set of points $X$ when $M$ moves along the line $BC$?}

\begin{figure}
\centering
\fbox{\includegraphics[width=8.5cm]{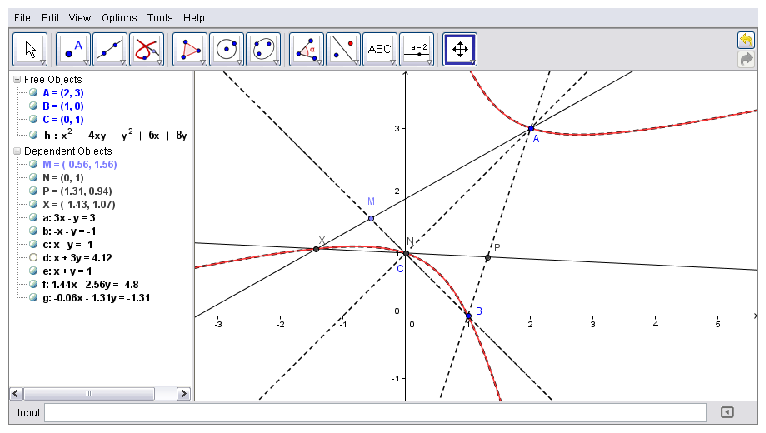}} \\
\fbox{\includegraphics[width=8.5cm]{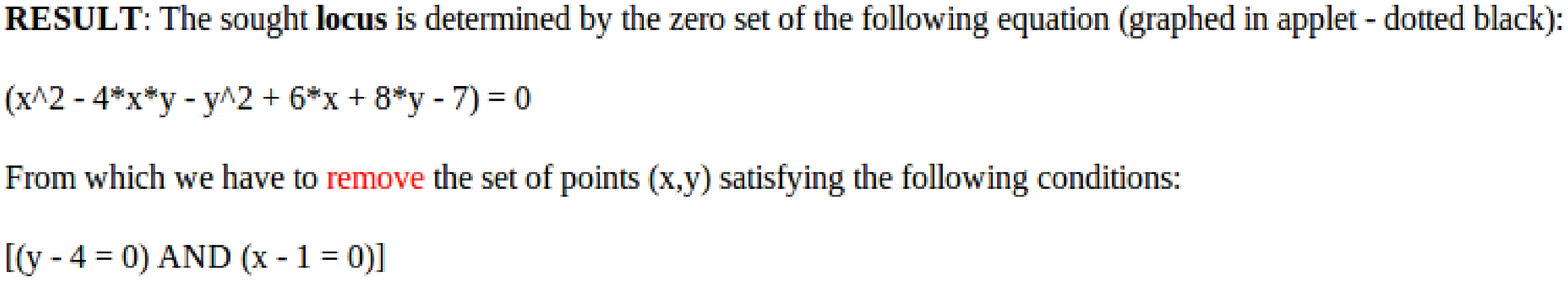}}
\caption{Locus as a constructible set.}
\label{fig:Locus-as-a-constructible-set-1}
\end{figure}

When the vertices of the triangle $ABC$ are the points $(2,3)$, $(1,0)$ and $(0,1)$, the locus set is a conic from which a point has to be removed. That is, the locus set is not an algebraic variety but a locally closed set. Figure \ref{fig:Locus-as-a-constructible-set-1} shows the plotting of the conic in GeoGebra together with its precise algebraic description as provided by the prototype.

If we consider this same construction with $A(0,0),B(1,0)$ and $C(0,1)$, a standard DGS will plot a straight line as locus, while ordinary elimination will give the true locus $2x+2y = 1$ plus two other lines, namely, the coordinate axes $x = 0$ and $y = 0$. These extra lines correspond to two degenerate positions for the mover: $M = B$ and $M = C$. Applying the criterion sketched above, the system identifies these two lines as special components and hence removes them from the final description. In tables \ref{table:outputCase1} and \ref{table:outputCase2} we find the parametric systems and outputs from the Gr\"obnerCover algorithm for the two considered instances respectively.

\begin{table}
\begin{center}
\begin{tabular}{|l|p{8cm} |}
   \hline 
   Parametric System & $-x_1-x_2+1, 2x_1+2x_2-2x_3-2x_4, -2x_3+2x_4-2, x_1+3x_2-x_5-3x_6, -3x_5+x_6+3, -(x_4-y)(x_3-x_5)+(x_4-x_6)(x_3-x), -(y-3)(x_1-2)+(x-2)(x_2-3)$
 \\
   \hline
   \hline
   Basis segment 1 & $\{1\}$\\
   \hline
   Segment 1 & $\mathbb{V}(0)\setminus \mathbb{V}(x^2-4xy+6x-y^2+8y-7)$  \\
   \hline
   Basis segment 2 $^*$ & $ \{\{(5y-20)x_6+(-3x+6y+3),(x-2y+7)x_6+(-3y)\},\{(5y-20)x_5+(-x-3y+21),(x-7y+27)x_5+(4y-28)\},x_4-1,x_3,\{(y-4)x_2+(-x+2y+1),(x-2y+7)x_2+(-5y)\},\{(y-4)x_1+(x-3y+3),(x-y+3)x_1+(4y-4)\}\}$ \\
   \hline
   Segment 2 & $\mathbb{V}(x^2-4xy+6x-y^2+8y-7)\setminus \mathbb{V}(y-4,x-1)$ \\
   \hline 
   Basis segment 3 & $\{1\}$ \\
   \hline
   Segment 3 & $\mathbb{V}(y-4,x-1)\setminus \mathbb{V}(1)$ \\
   \hline
   \hline
   Locus (after heuristic step) & $\mathbb{V}(x^2 - 4xy - y^2 + 6x + 8y - 7) \setminus \mathbb{V}(y - 4,x - 1)$ \\
   \hline
 \end{tabular}
\smallskip
\smallskip
\caption{Parametric system, Gr\"obnerCover output and returned locus for instance with $A(2,3)$, $B(1,0)$ and $C(0,1)$. \newline \footnotesize{$^*$ Includes regular functions (see \cite{MontesWibmer2010}).}}
\label{table:outputCase1}
\end{center}
\end{table}

\begin{table}
\begin{center}
\begin{tabular}{|l|p{8cm}|}
   \hline 
   Parametric System & $-x_1-x_2+1,x_3, -x_2 + x_4, -x_1 + x_5, -x_6, -x_1 y + x_2 x, -(x_4 - y)(x_3 - x_5) + (x_4 - x_6)(x_3 - x)$ \\
   \hline
   \hline
   Basis segment 1 & $\{1\}$ \\
   \hline
   Segment 1 & $\mathbb{V}(0)\setminus (\mathbb{V}(2x+2y-1) \cup \mathbb{V}(x) \cup \mathbb{V}(y))$ \\
   \hline
   Basis segment 2 & $\{x_6,(x+y)x_5-x,(x+y)x_4-y,x_3,(x+y)x_2-y,(x+y)x_1-x\}$ \\
   \hline
   Segment 2 & $(\mathbb{V}(2x+2y-1)\setminus \mathbb{V}(1)) \cup (\mathbb{V}(x)\setminus \mathbb{V}(x,y)) \cup (\mathbb{V}(y)\setminus \mathbb{V}(x,y))$ \\
   \hline
   Basis segment 3 & $\{x_6,x_4+x_5-1,x_3,x_2+x_5-1,x_1-x_5,x_5^2-x_5\}$ \\
   \hline
   Segment 3 & $ \mathbb{V}(x,y)\setminus \mathbb{V}(1)$ \\
   \hline
   \hline
   Locus (after heuristic step) & $\mathbb{V}(2x + 2y - 1)$ \\
   \hline  
 \end{tabular}
\smallskip
\smallskip
\caption{Parametric system, Gr\"obnerCover output and returned locus for instance with $A(0,0)$, $B(1,0)$ and $C(0,1)$.}
\label{table:outputCase2}
\end{center}
\end{table}

\section{Prototype description}

The system (accessible at \cite{LocusEnvelopeWeb}) consists of a web page with a \mbox{GeoGebra} applet where the user constructs a locus or a family of linear objects depending on a point. For any of these constructions (specified using a predetermined set of GeoGebra commands) the prototype provides the algebraic description of the locus/envelope by just pressing one button. Note that in its current state, the system does not provide the equation of the envelope, but the one of the discriminant line. A note stating this fact should be given if using the system for teaching purposes. 

The process is roughly as follows. First, the XML description of the GeoGebra construction is sent to a Server where an installation of a Sage Cell Server (\cite{SageCellServer}) is maintained by the authors. There, the construction follows an algebraization process, as specified by a Sage library \cite{BotanaICCSA2011}. The communication Sage-GeoGebra is made possible by the JavaScript GeoGebra functions that allow the data transmission to and from the applet. In particular, the XML description of any GeoGebra diagram can be obtained. The processing of the XML description of the diagram is made by some ad-hoc code by the authors that use Sage through the Sage cell server, a general service by Sage. More concretely, Singular, included in Sage and with an implementation of the Gr\"obnerCover algorithm, is used. The Gr\"obner cover of the obtained parametric polynomial system is analyzed, and the accepted components of the locus/envelope are incorporated into the applet. 



Note that the goal is not to provide a final tool but a proof-of-concept prototype showing the feasibility of using sophisticated algorithms like Gr\"obnerCover to supplement the symbolic capabilities of existing dynamic geometry systems, as well as to show the advantage of connecting different systems by using web services.

\bibliographystyle{splncs} 
\bibliography{AutomaticDeduction}

\end{document}